\documentclass[english,reqno]{amsart}
\usepackage{pslatex}
\usepackage[T1]{fontenc}
\usepackage{graphicx}
\usepackage{amssymb}

\makeatletter
 \theoremstyle{plain}    
 \newtheorem{thm}{Theorem}[section]
 \numberwithin{equation}{section} 
 \numberwithin{figure}{section} 
 \theoremstyle{plain}
 \theoremstyle{plain}    
 \newtheorem{prop}[thm]{Proposition} 
 \theoremstyle{plain}    
 \newtheorem{cor}[thm]{Corollary} 
 \theoremstyle{plain}    
 \newtheorem{lem}[thm]{Lemma} 

\usepackage{psfrag}
\usepackage[colorlinks=false,bookmarks,pdfhighlight=/O,pdfauthor={Marc Kesseböhmer},pdftitle={IES}]{hyperref}

\DeclareMathOperator{\ind}{\boldsymbol{1}}

\newcommand{\1}{\boldsymbol{1}}

\newcommand{\ep}{\hfill $\Box$} 
\newcommand{\proc}[1]{\medbreak \noindent{\em{#1}}} 

\AtBeginDocument{
  
}

\usepackage{babel}
\makeatother
\begin{document}
\title[Limit laws for distorted return time processes]{Limit laws for distorted return time processes for infinite measure preserving transformations}

\author{Marc Kesseb{\"O}hmer and Mehdi Slassi}

\email{mhk@math.uni-bremen.de, slassi@math.uni-bremen.de}

\address{Universit{\"a}t Bremen, Fachbereich 3 f{\"u}r Mathematik und Informatik, Bibliothekstra{\ss}e~1,\ D--28359 Bremen, Germany.}


\date{23~March~2006}

\thanks{Research supported by Zentrale Forschungsförderung Universit{\"a}t Bremen}


\begin{abstract} We consider conservative ergodic measure preserving transformations on infinite measure spaces and investigate the asymptotic behaviour of  distorted return time processes with respect to sets  satisfying a type of Darling-Kac condition. We identify two critical cases for which we prove uniform distribution laws. For this we introduce the notion of uniformly returning sets and discuss some of their properties.  

\end{abstract}

\maketitle
\let\languagename\relax

\section{Introduction and statement of main results}

In this paper $\left(X,T,\mathcal{A},\mu\right)$ will always denote
a conservative ergodic measure preserving dynamical systems where
$\mu$ is an infinite $\sigma$-finite measure. In particular, this
means that the mean return time to sets of finite positive measure
is infinite. Hence if the system is given by a Markov chain, this
corresponds to the \emph{null recurrent} situation. The investigation
of ergodic and probabilistic properties of such dynamical systems
leads to a number of interesting results which can always be interpreted
within the theory of null recurrent Markov chains and which sometimes
generalize classical theorems within this theory.

In this paper we present a generalization of the Thaler-Dynkin-Lamperti
arc-sine law (cf. \cite{Thaler:98} and (T) in Subsection \ref{sub:Limit-laws.})
describing the asymptotic behaviour of the renewal theoretic process
$Z_{n}$ given by \[
Z_{n}(x):=\left\{ \begin{array}{ll}
\max\{ k\leq n:\; T^{k}(x)\in A\}, & x\in A_{n}:=\bigcup_{k=0}^{n}T^{-k}A,\\
0, & \textrm{else.}\end{array}\right.\]
For a regularly varying function $F$ we consider the \emph{distorted
processes} \[
\frac{F\left(Z_{n}\right)}{F\left(n\right)}\quad\textrm{and }\quad\frac{F\left(n-Z_{n}\right)}{F\left(n\right)}.\]
 In particular we introduce the processes\[
\Phi_{n}:=\frac{\sum_{k=0}^{Z_{n}}\mu\left(A\cap\left\{ \varphi>k\right\} \right)}{\mu\left(A_{n}\right)}\quad\textrm{and}\quad\Psi_{n}:=\frac{\sum_{k=0}^{n-Z_{n}}\mu\left(A\cap\left\{ \varphi>k\right\} \right)}{\mu\left(A_{n}\right)},\]
 which we refer to as the \emph{normalized Kac process} and \emph{normalized
spent time Kac process,} respectively. In here, \begin{equation}
\varphi(x)=\inf\{ n\geq1:\; T^{n}(x)\in A\},\;\; x\in X,\label{returntime}\end{equation}
 denotes the \emph{first return time} to the set $A$.

In Proposition \ref{propo0} we give a purely probability theoretical
result allowing us to derive limit laws of distorted processes if
the limit law for the corresponding original process is known. This
result is then applied in Corollary \ref{coro 1} to treat $\Phi_{n}$
and $\Psi_{n}$. 

Two critical cases are identified which are not covered by Proposition
\ref{propo0} and are subject of our two main Theorems \ref{theo1}
and \ref{theo2}. More precisely this means that if the sequence $\left(\frac{Y_{n}}{n}\right)$
converges in distribution to $0$ and $L$ is slowly varying, then
in general nothing is known about the asymptotics of $\frac{L\left(Y_{n}\right)}{L\left(n\right)}$.
In this situation we are able to show that under suitable conditions
on the wandering rate of a uniform set $A$ we have\[
\frac{L\left(Z_{n}\right)}{L\left(n\right)}\stackrel{\mathcal{L\left(\mu\right)}}{\;\Longrightarrow}\;\boldsymbol{U},\]
and if $A$ is a uniformly returning set -- as introduced in Subsection
\ref{sub:Uniform-and-uniform} -- we have\[
\Psi_{n}\stackrel{\mathcal{L\left(\mu\right)}}{\;\Longrightarrow}\;\boldsymbol{U}.\]
In here, $\boldsymbol{U}$ denotes a random variable \emph{distributed
uniformly} on $\left[0,1\right]$. 

Obviously, Theorem \ref{theo1} can be applied to infinite measure
preserving interval maps $T:\left[0,1\right]\longrightarrow\left[0,1\right]$
with indifferent fixed points satisfying the Thaler condition stated
in \cite{thaler:95}, whereas Theorem \ref{theo2} is applicable to
those map satisfying the corresponding condition in \cite{Thaler:00}.
Other examples in the context of continued fractions -- also covered
by Theorem \ref{theo2} -- are treated in \cite{KessSlassi:05}. For
related results we refer to \cite{ThalerZweimuell:04} and for further
interesting results concerning distributional limit theorems for ergodic
sums in this context to \cite{Zweimueller:03}.

\subsection{Infinite ergodic theory. \label{sub:Infinite-ergodic-theory} }

A characterization of $\left(X,T,\mathcal{A},\mu\right)$ being a
con\-ser\-vative ergodic measure preserving dynamical system where
$\mu$ is an infinite $\sigma$-finite measure as used in this paper
will be given at the end of this subsection. For further definitions
and details we refer the reader to \cite{Aaronson:97}.

Let \[
\mathcal{P}_{\mu}:=\left\{ \nu:\nu\:\textrm{probability measure on}\:\mathcal{A}\,\textrm{with }\nu\ll\mu\right\} \]
 denote the set of probability measures on $\mathcal{A}$ which are
absolutely continuous with respect to $\mu$. The measures from $\mathcal{P}_{\mu}$
represent the admissible initial distributions for the processes associated
with the iteration of $T$. The symbol $\mathcal{P}_{\mu}$ will also
be used for the set of the corresponding densities. 

Let us recall the notion of the wandering rate. For a fixed set $A\in\mathcal{A}$
with $0<\mu\left(A\right)<\infty$ we set 

\[
A_{n}:=\bigcup_{k=0}^{n}T^{-k}A\quad\mathrm{and}\quad W_{n}:=W_{n}\left(A\right):=\mu\left(A_{n}\right),\qquad n\geq0,\]
and call the sequence $\left(W_{n}\left(A\right)\right)$ the \emph{wandering
rate} of $A.$ Note that for the wandering rate the following identities
hold\[
W_{n}\left(A\right)=\sum_{k=0}^{n}\mu\left(A\cap\{\varphi>k\}\right)=\int_{A}\min(\varphi,n+1)\, d\mu\\
,\qquad n\geq0.\]
 Since $T$ is conservative and ergodic, for all $\nu\in\mathcal{P_{\mu}}$,\[
\lim_{n\to\infty}\nu\left(A_{n}\right)=1\quad\textrm{and}\quad\nu\left(\left\{ \varphi<\infty\right\} \right)=1.\]

The key to an understanding of the stochastic properties of a nonsingular
transformation of a $\sigma$--finite measure space often lies in
the study of the long-term behaviour of the iterates of its \emph{transfer
operator} \[
\hat{T}:L_{1}\left(\mu\right)\longrightarrow L_{1}\left(\mu\right),\; f\longmapsto\hat{T}\left(f\right):=\frac{d\left(\nu_{f}\circ T^{-1}\right)}{d\mu},\]
where $\nu_{f}$ denote the measure with density $f$ with respect
to $\mu$. Clearly, $\hat{T}$ is a positive linear operator characterized
by\[
\int_{B}\hat{T}\left(f\right)\; d\mu=\int_{T^{-1}\left(B\right)}f\; d\mu,\qquad f\in L_{1}\left(\mu\right),\quad B\in\mathcal{A}.\]
An approximation argument shows that equivalently for all $f\in L_{1}\left(\mu\right)$
and $g\in L_{\infty}\left(\mu\right)$\[
\int_{X}\hat{T}\left(f\right)\cdot g\; d\mu=\int_{X}f\cdot g\circ T\; d\mu.\]

The ergodic properties of $(X,T,\mathcal{A},\mu)$ can be characterized
in terms of the transfer operator in the following way (cf. \cite[Proposition 1.3.2]{Aaronson:97}).
A system is conservative and ergodic if and only if for all $f\in L_{1}^{+}\left(\mu\right):=\left\{ f\in L_{1}\left(\mu\right):\; f\geq0\;\mathrm{and}\;\int_{X}f\; d\mu>0\right\} $
we have $\mu$-a.e. \begin{equation}
\sum_{n\geq0}\hat{T}^{n}\left(f\right)=\infty.\label{eq:conserErgod}\end{equation}
 Invariance of $\mu$ under $T$ means $\hat{T}\left(1\right)=1.$

\subsection{Uniform and uniformly returning sets\label{sub:Uniform-and-uniform}}

The following two definitions are in many situation crucial within
infinite ergodic theory.

\begin{itemize}
\item A set $A\in\mathcal{A}$ with $0<\mu\left(A\right)<\infty$ is called
\emph{uniform for} $f\in\mathcal{P_{\mu}}$ if there exists a sequence
$\left(a_{n}\right)$ of positive reals such that\[
\frac{1}{a_{n}}\sum_{k=0}^{n-1}\hat{T}^{k}\left(f\right)\;\longrightarrow\;1\qquad\mu-\textrm{a.e. \; uniformly\; on}\; A\]
(i.e. uniform convergence in $L_{\infty}\left(\mu|_{A\cap\mathcal{A}}\right)$).
\item The set $A$ is called a \emph{uniform} set \emph{}if it is uniform
for some $f\in\mathcal{P_{\mu}}.$
\end{itemize}
\proc{Remark.} \label{remark0}Note that from \cite[Proposition 3.8.7]{Aaronson:97}
we know, that $\left(b_{n}\right)$ is regularly varying with exponent
$\alpha$ (for the definition of this property see Section \ref{sub:Classical-results-on})
if and only if $\left(W_{n}\right)$ is regularly varying with exponent
$\left(1-\alpha\right)$. In this case $\alpha$ lies in the interval
$\left[0,1\right]$ and \begin{equation}
a_{n}W_{n}\sim\frac{n}{\Gamma\left(1+\alpha\right)\Gamma\left(2-\alpha\right)}.\label{eq:AsymptAaronson}\end{equation}
In here, $c_{n}\sim a_{n}$ for some sequences $\left(c_{n}\right)$
and $\left(a_{n}\right)$ means that $a_{n}\not=0$ has only finitely
many exceptions and $\lim_{n\to\infty}\frac{c_{n}}{a_{n}}=1$.

Next, we define a new property for sets similar to that of being uniform.
It will be used to state the conditions in Theorem \ref{theo2}.

\proc{Definition.} A set $A\in\mathcal{A}$ with $0<\mu\left(A\right)<\infty$
is called \emph{uniformly returning} \emph{for} $f\in\mathcal{P_{\mu}}$
if there exists a positive increasing sequence $\left(b_{n}\right)$
diverging to $\infty$ such that \[
b_{n}\hat{T}^{n}\left(f\right)\;\longrightarrow\;1\quad\mu-\textrm{a.e. \; uniformly\; on}\; A.\]
The set $A$ is \emph{}called \emph{uniformly} \emph{returning} if
it is uniformly returning for some $f\in\mathcal{P_{\mu}}.$\medbreak

The following Example show the existence of uniformly returning sets.

\proc{Example.}  Let $T:\left[0,1\right]\longrightarrow\left[0,1\right]$
be an interval map with two increasing full branches and an indifferent
fixed point at $0$ satisfying Thaler's conditions in \cite{Thaler:00}.
Then any set $A\in\mathcal{B}_{\left[0,1\right]}$ with positive distance
from the indifferent fixed point $0$ and $\lambda\left(A\right)>0$
is uniformly returning. As an typical example see the Lasota--York
interval map (cf. Example after Theorem \ref{theo2}).

\begin{prop}
\label{pro:URimplU}Any uniformly returning set is uniform. 
\end{prop}
\proc{Proof.} Let $A$ be a uniformly returning set for $f\in\mathcal{P}_{\mu}$.
Then for each \emph{$\varepsilon\in\left(0,1\right)$ there exists}
a positive integer $n_{0}$ such that for all $n\geq n_{0}$ we have\[
\left(1-\varepsilon\right)\frac{1}{b_{n}}\leq\hat{T}^{n}\left(f\right)\leq\left(1+\varepsilon\right)\frac{1}{b_{n}}\quad\mu\textrm{-a.e.\; uniformly on }\; A.\]
 Set $\widetilde{f}:=\hat{T}^{n_{0}}\left(f\right)\in\mathcal{P}_{\mu}$.
Since $b_{n}\uparrow\infty$ we deduce that $\left(\hat{T}^{n}\left(\widetilde{f}\right)\right)_{n\geq0}$
is uniformly bounded $\mu$-a.e. on $A$. Furthermore, since $T$
is conservative and ergodic we have by (\ref{eq:conserErgod}) \[
\sum_{k=0}^{\infty}\hat{T}^{k}\left(\tilde{f}\right)=\infty\qquad\mu\textrm{-a.e.}\]
Using the fact that \[
b_{n+n_{0}}\hat{T}^{n}\left(\widetilde{f}\right)\;\longrightarrow\;1\quad\mu\textrm{-a.e.\; uniformly on }\; A\]
and that $\left(\hat{T}^{n}\left(\widetilde{f}\right)\right)_{n\geq0}$
is uniformly bounded, we get \[
\frac{1}{\sum_{k=0}^{n}\frac{1}{b_{k}}}\cdot\sum_{k=0}^{n}\hat{T}^{k}\left(\widetilde{f}\right)\;\longrightarrow\;1\quad\mu\textrm{-a.e.\; uniformly on }\; A.\]
This shows that $A$ is a uniform set for $\widetilde{f}$.\ep\medbreak

\proc{Remark.} The inverse implication of Proposition \ref{pro:URimplU}
is stated in \cite{KessSlassi:05} under some additional assumptions. 

To characterize the difference of the notions of uniform and uniformly
returning set we make the following considerations. Suppose $A$ is
uniformly returning for $f:=\frac{1}{\mu\left(B\right)}\ind_{B}$,
$B\in\mathcal{A}$, $0<\mu\left(B\right)<\infty$, then there exists
a sequence $\left(b_{n}\right)$ such that \[
\frac{b_{n}}{\mu\left(B\right)}\hat{T}^{n}\left(\ind_{B}\right)\to1\qquad\mu\textrm{-a.e.\; uniformly on }\; A.\]
Integrating over $A$ yields\[
b_{n}\mu\left(B\cap T^{-n}A\right)\to\mu\left(B\right)\mu\left(A\right).\]
In an analog way we deduce for a uniform set $A$ that \[
\frac{1}{a_{n}}\sum_{k=1}^{n}\mu\left(B\cap T^{-k}A\right)\to\mu\left(B\right)\mu\left(A\right).\]
Hence, we may interpret '\emph{uniform returning}' as a version of
'\emph{strong mixing}' in infinite ergodic theory, whereas \emph{'uniform'}
corresponds to a version of '\emph{ergodicity}'. 

Next we characterize the sequence $\left(b_{n}\right)$ in this definition
by the wandering rate similarly to (\ref{eq:AsymptAaronson}).

\begin{prop}
\label{neue pro} For $\beta\in\left[0,1\right)$ we have that $\left(b_{n}\right)$
is regularly varying with exponent $\beta$ if and only if $\left(W_{n}\right)$
is regularly varying with the same exponent. In this case,\[
b_{n}\sim W_{n}\Gamma\left(1-\beta\right)\Gamma\left(1+\beta\right)\qquad\left(n\to\infty\right).\]

\end{prop}
The proof of this proposition and of the following proposition and
theorems will be postponed to Section \ref{sec:Proof-of-main}.

\subsection{Limit laws.\label{sub:Limit-laws.}}

An important question when studying convergence in distribution for
processes defined in terms of a non-singular transformation is to
what extent the limiting behaviour depends on the initial distribution.
This is formalized as follows.

Let $\nu$ be a probability measure on the measurable space $(X,\mathcal{A})$
and $\left(R_{n}\right)_{n\geq1}$ be a sequence of measurable real
functions on $X$, distributional convergence of $\left(R_{n}\right)_{n\geq1}$
w.r.t. $\nu$ to some random variable $R$ with values in $\left[-\infty,\infty\right]$
will be denoted by $R_{n}\stackrel{\nu}{\Longrightarrow}R$. \emph{Strong
distributional convergence} abbreviated by $R_{n}\stackrel{\mathcal{L\left(\mu\right)}}{\Longrightarrow}R$
on the $\sigma$--finite measures space $\left(X,\mathcal{A,\mu}\right)$
means that $R_{n}\stackrel{\nu}{\Longrightarrow}R$ for all $\nu\in\mathcal{P_{\mu}}$.
In particular for $c\in\left[-\infty,\infty\right]$,\[
R_{n}\stackrel{\mathcal{L\left(\mu\right)}}{\Longrightarrow}c\quad\Longleftrightarrow\quad R_{n}\longrightarrow c\textrm{\quad locally in measure,}\]
which we also denote by $R_{n}\stackrel{\mu}{\longrightarrow}c$.

Now we are in the position to state the first interesting limit law
for the process $Z_{n}$ which is due to Thaler \cite{Thaler:98}.

\begin{itemize}
\item [(T)] \textbf{Thaler's Dynkin-Lamperti arc-sin Law.} Let $A\in\mathcal{A}$
with $0<\mu\left(A\right)<\infty$ be a uniform set. If the wandering
rate $\left(W_{n}\left(A\right)\right)$ is regularly varying with
exponent $1-\alpha$ for $\alpha\in\left[0,1\right]$, then we have\begin{equation}
\frac{Z_{n}}{n}\;\stackrel{\mathcal{L\left(\mu\right)}}{\Longrightarrow}\;\xi_{\alpha}.\label{tharesult}\end{equation}
In here, for $\alpha\in\left(0,1\right),$ $\xi_{\alpha}$ denotes
the random variable on $\left[0,1\right]$ with density \[
f_{\xi_{\alpha}}\left(x\right)=\frac{\sin\pi\alpha}{\pi}\frac{1}{x^{1-\alpha}\left(1-x\right)^{\alpha}},\qquad0<x<1.\]
The distribution of $\xi_{\alpha}$ is also called \emph{the generalized
arc-sine distribution.} The continuous extension is given by \emph{$\xi_{0}=0$}
and \emph{$\xi_{1}=1.$}
\end{itemize}
To apply (T) to the distorted processes we need the following proposition.

\begin{prop}
\label{propo0} Let $\left(\Omega,\mathcal{F},\mathbb{P}\right)$
be a probability space, let $Y_{n}:\Omega\longrightarrow\left[0,\infty\right]$
be measurable $\left(n\geq1\right)$, and let $Y$ be a random variable
with values in $\left[0,\infty\right]$.
\begin{itemize}
\item [{\rm (1)}]If $\mathbb{P}\left(Y=0\right)=0=\mathbb{P}\left(Y=\infty\right)$
and $F$ is a regularly varying function with exponent $\beta\in\mathbb{R},$
then \[
\frac{Y_{n}}{n}\stackrel{\mathbb{P}}{\Longrightarrow}Y\qquad\quad\implies\qquad\quad\frac{F\left(Y_{n}\right)}{F\left(n\right)}\stackrel{\mathbb{P}}{\Longrightarrow}Y^{\beta}.\]

\item [{\rm (2)}]If $Y=0$ and $F$ is a regularly varying function with
exponent $\beta\in\mathbb{R}\setminus\left\{ 0\right\} $ then\[
\frac{Y_{n}}{n}\stackrel{\mathbb{P}}{\Longrightarrow}0\qquad\quad\implies\qquad\quad\frac{F\left(Y_{n}\right)}{F\left(n\right)}\stackrel{\mathbb{P}}{\Longrightarrow}\left\{ \begin{array}{ll}
0 & \textrm{for }\beta>0\\
\infty & \textrm{for }\beta<0\end{array}\right..\]
 
\end{itemize}
\end{prop}
\begin{figure}
\psfrag{00.01}{{\small $\alpha=0.01$}}\psfrag{00.3}{{\small $\alpha=0.3$}}\psfrag{00.5}{{\small $\alpha=0.5$}}\psfrag{00.6}{{\small $\alpha=0.6$}}\psfrag{00.7}{{\small $\alpha=0.7$}}\psfrag{00.9}{{\small $\alpha=0.9$}}\psfrag{02}{{\small $2$}}\psfrag{00}{{\small $0$}}\psfrag{01}{{\small $1$}}

\includegraphics[%
  width=1.0\textwidth]{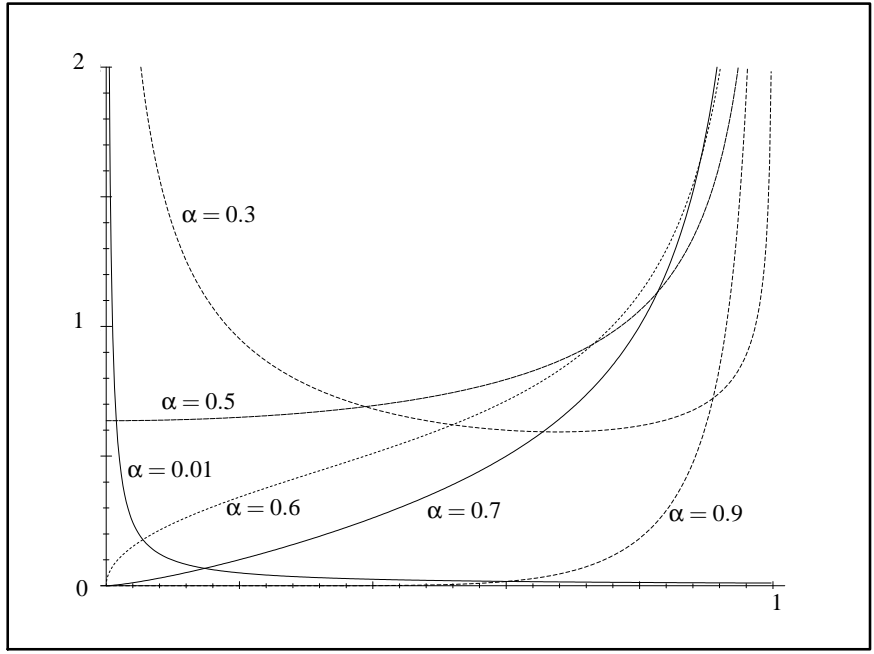}

\caption{\label{cap:Figdensity1} The densities $f_{X_{\alpha}}$ of the limiting
distribution of the normalized Kac process for different values of
$\alpha\in\left(0,1\right)$. The extreme distributions corresponding
to $\alpha=0$ and $1$ are the Dirac measures $\delta_{0}$ and $\delta_{1}$,
respectively.}
\end{figure}
The following corollary is a direct consequence of (\ref{tharesult}),
Proposition\ref{propo0}, and the fact that $\Phi_{n}=\frac{F\left(Z_{n}\right)}{F\left(n\right)},\;\Psi_{n}=\frac{F\left(n-Z_{n}\right)}{F\left(n\right)}$
with $F\left(n\right):=W_{n}$. 

\begin{figure}
\psfrag{00.01}{{\small $\alpha=0.01$}}\psfrag{00.3}{{\small $\alpha=0.3$}}\psfrag{00.5}{{\small $\alpha=0.5$}}\psfrag{00.2}{{\small $\alpha=0.2$}}\psfrag{00.7}{{\small $\alpha=0.7$}}\psfrag{00.97}{{\small $\alpha=0.97$}}\psfrag{02}{{\small $2$}}\psfrag{00}{{\small $0$}}\psfrag{01}{{\small $1$}}

\includegraphics[%
  width=1.0\textwidth]{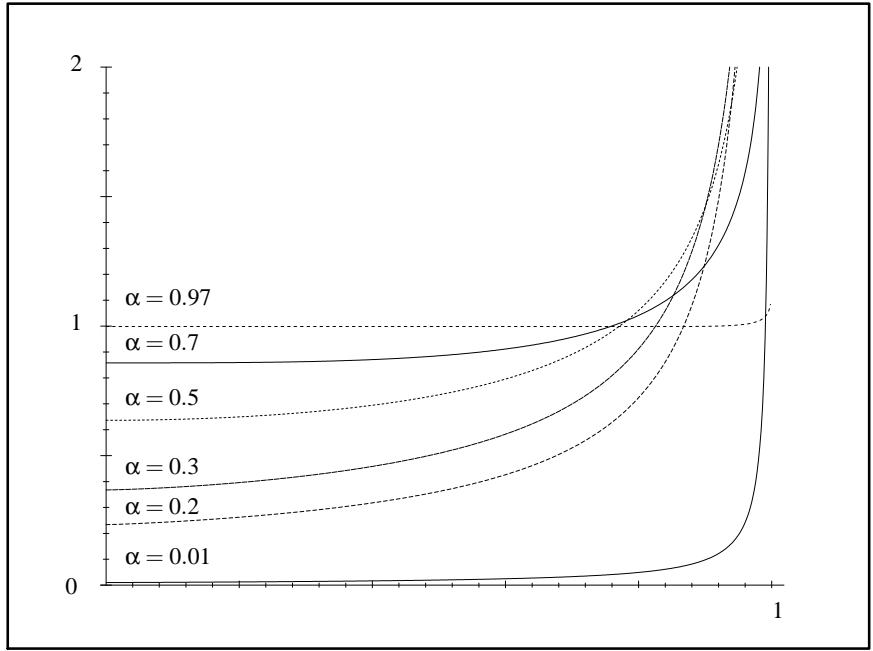}

\caption{\label{cap:Figdensity2} The densities $f_{Y_{\alpha}}$ of the limiting
distribution of the normalized spent time Kac process for different
values of $\alpha\in\left(0,1\right)$. The extreme distribution corresponding
to $\alpha=0$ and $1$ are the Dirac measure $\delta_{1}$ and the
uniform distribution on $\left[0,1\right]$, respectively.}
\end{figure}

\begin{cor}
\label{coro 1}Let $A\in\mathcal{A}$ with $0<\mu\left(A\right)<\infty$
be a uniform set. If the wandering rate $\left(W_{n}\right)$ is regularly
varying with exponent $1-\alpha$, then 
\begin{itemize}
\item [{\rm (1)}] If $0\leq\alpha\leq1,$ then we have\[
\Phi_{n}\;\stackrel{\mathcal{L\left(\mu\right)}}{\Longrightarrow}\; X_{\alpha},\]
where $X_{\alpha}$ denotes the random variable on $\left[0,1\right]$
with density \[
f_{X_{\alpha}}\left(x\right)=\frac{1}{1-\alpha}\frac{\sin\pi\alpha}{\pi}\frac{1}{x^{\frac{1-2\alpha}{1-\alpha}}\left(1-x^{\frac{1}{1-\alpha}}\right)^{\alpha}},\quad\alpha\in\left(0,1\right)\]
and $X_{0}=0$, $X_{1}=1$ (cf. Fig. \ref{cap:Figdensity1}). 
\item [{\rm (2)}] If $0\leq\alpha<1,$ then we have\[
\Psi_{n}\;\stackrel{\mathcal{L\left(\mu\right)}}{\Longrightarrow}\; Y_{\alpha},\]
where $Y_{\alpha}$ denotes the random variable on $\left[0,1\right]$
with density \[
f_{Y_{\alpha}}\left(x\right)=\frac{1}{1-\alpha}\frac{\sin\pi\alpha}{\pi}\frac{1}{\left(1-x^{\frac{1}{1-\alpha}}\right)^{1-\alpha}},\quad\alpha\in\left(0,1\right)\]
and $Y_{0}=1$ (cf. Fig. \ref{cap:Figdensity2}).
\end{itemize}
\end{cor}
\proc{Remark.}  For $\alpha\in\left(0,1\right)$ we have \[
X_{\alpha}\stackrel{\textrm{dist.}}{=}\left(\xi_{\alpha}\right)^{1-\alpha}\quad\textrm{and}\qquad Y_{\alpha}\stackrel{\textrm{dist.}}{=}\left(1-\xi_{\alpha}\right)^{1-\alpha}.\]
Note, that in particular both $X_{\frac{1}{2}}$ and $Y_{\frac{1}{2}}$obey
the arc-sine law, i.e. they have density\[
f_{Y_{\frac{1}{2}}}\left(x\right)=f_{X_{\frac{1}{2}}}\left(x\right)=\frac{2}{\pi}\frac{1}{\left(1-x^{2}\right)^{\frac{1}{2}}},\qquad0<x<1.\]

Now we state as our first main result the following \emph{uniform
distribution law}. Note that this corresponds to the case $\beta=0$
and $Y=0$ in Proposition \ref{propo0}, which therefore is not applicable. 

\begin{thm}
\label{theo1}Let $A\in\mathcal{A}$ with $0<\mu\left(A\right)<\infty$
be a uniform set. If the wandering rate $\left(W_{n}\right)$ is regularly
varying with exponent $1$ such that $W_{n}\sim n\cdot\frac{1}{L\left(n\right)}$,
where $L$ is a slowly varying function . Then we have\[
\frac{L\left(Z_{n}\right)}{L\left(n\right)}\stackrel{\mathcal{L\left(\mu\right)}}{\;\longrightarrow}\;\boldsymbol{U},\]
where $\boldsymbol{U}$ denotes the uniformly distributed random variable
on $\left[0,1\right]$.
\end{thm}
\proc{Remark.} Under the assumptions of Theorem \ref{theo1} we have
by (T) and the first part of Corollary \ref{coro 1} that\[
\frac{Z_{n}}{n}\stackrel{\mathcal{L\left(\mu\right)}}{\longrightarrow}0\quad\textrm{and}\quad\frac{1}{\Phi_{n}}\stackrel{\mathcal{L\left(\mu\right)}}{\longrightarrow}\infty.\]
However, for the product of the two processes we have by Theorem \ref{theo1}\[
\frac{Z_{n}}{n\Phi_{n}}\stackrel{\mathcal{L\left(\mu\right)}}{\;\longrightarrow}\;\boldsymbol{U}.\]

\proc{Example.} Let $f\left(0\right)=0,\; f\left(x\right)=x+x^{2}e^{-\frac{1}{x}},\quad x>0,$
and let $a\in\left(0,1\right)$ be determined by $f\left(a\right)=1.$
Define $T:\left[0,1\right]\longrightarrow\left[0,1\right]$ by\[
T\left(x\right):=\left\{ \begin{array}{ll}
f\left(x\right), & x\in\left[0,a\right],\\
\frac{x-a}{1-a}, & x\in\left(a,1\right].\end{array}\right.\]
Then the map $T$ satisfies Thaler's conditions (T1)--(T4) in \cite{thaler:95}.
Any set $A\in\mathcal{B}_{\left[0,1\right]}$ with $\lambda\left(A\right)>0$
which is bounded away from the indifferent fixed points is a uniform
set for $T.$ Furthermore we have \[
W_{n}\sim const\cdot\frac{n}{\log\left(n\right)}\qquad\left(n\to\infty\right).\]
 Hence,\[
\frac{\log Z_{n}}{\log n}\:\stackrel{\mathcal{L\left(\mu\right)}}{\Longrightarrow}\;\boldsymbol{U}.\]

Next we state as our second main result the uniform distribution law
for the normalized spent time Kac process. This again corresponds
to the case $\beta=0$ and $Y=0$ in Proposition \ref{propo0}. 

\begin{thm}
\label{theo2} Let $A\in\mathcal{A}$ with $0<\mu\left(A\right)<\infty$
be a uniformly returning set$.$ If the wandering rate $\left(W_{n}\right)$
is slowly varying, then we have\[
\Psi_{n}\stackrel{\mathcal{L\left(\mu\right)}}{\;\Longrightarrow}\;\boldsymbol{U},\]
where the random variable $\boldsymbol{U}$ is distributed uniformly
on $\left[0,1\right]$.
\end{thm}
\proc{Example.}  We consider the \emph{Lasota--Yorke} map $T:\left[0,1\right]\longrightarrow\left[0,1\right]$,
defined by\[
T\left(x\right):=\left\{ \begin{array}{ll}
\frac{x}{1-x}, & x\in\left[0,\frac{1}{2}\right],\\
2x-1, & x\in\left(\frac{1}{2},1\right].\end{array}\right.\]
This map satisfies the Thaler's conditions (i)--(iv) in \cite{Thaler:00}.
Any compact subset $A$ of $\left(0,1\right]$ with $\lambda\left(A\right)>0$
is a uniformly returning set and we have\[
W_{n}\sim\log\left(n\right)\quad\textrm{as}\quad n\to\infty.\]
Hence,\[
\frac{\log\left(n-Z_{n}\right)}{\log\left(n\right)}\stackrel{\mathcal{L\left(\mu\right)}}{\;\Longrightarrow}\;\boldsymbol{U}.\]

\section{Regular variation and Tauberian results\label{sub:Classical-results-on}}

In this section we give some preparatory facts and results needed
for the proofs of the main statements in the following Section \ref{sec:Proof-of-main}. 

We first recall the concepts of regularly varying functions and sequences
(see also \cite{BinghamGoldieTeugels:89} for a comprehensive account).
Throughout we use the convention that for two sequences $\left(a_{n}\right)$,
$\left(b_{n}\right)$ we write $a_{n}=o\left(b_{n}\right)$ if $b_{n}\not=0$
fails only for finitely many $n$ $\textrm{and}\;\lim_{n\to\infty}\left|\frac{a_{n}}{b_{n}}\right|=0$.

A measurable function $R:\mathbb{R}^{+}\rightarrow\mathbb{R}$ with
$R>0$ on $\left(a,\infty\right)$ for some $a>0$ is called \emph{regularly
varying} at $\infty$ with exponent $\rho\in\mathbb{R}$ if\[
\lim_{t\to\infty}\frac{R\left(\lambda t\right)}{R\left(t\right)}=\lambda^{\rho}\quad{\rm for\; all}\;\lambda>0.\]
A regularly varying function $L$ with exponent $\rho=0$ is called
\emph{slowly varying} at $\infty$, i.e.\[
\lim_{t\to\infty}\frac{L\left(\lambda t\right)}{L\left(t\right)}=1\quad{\rm for\; all}\;\lambda>0.\]
 Clearly, a function $R:\mathbb{R}^{+}\rightarrow\mathbb{R}$ is  regularly
varying \emph{at} $\infty$ with exponent $\rho$$\in\mathbb{R}$
if and only if\[
R\left(t\right)=t^{\rho}L\left(t\right),\quad t\in\mathbb{R}^{+},\]
for $L$ slowly varying at $\infty.$

A function $R$ is said to be \emph{regularly varying at $0$} if
$t\mapsto R\left(\frac{1}{t}\right)$ is regularly varying at $\infty.$

A \emph{sequence} $\left(u_{n}\right)$ is \emph{regularly varying
with exponent} $\rho$ if $u_{n}=R\left(n\right),\; n\geq1,$ for
$R:\mathbb{R}^{+}\rightarrow\mathbb{R}$ regularly varying at $\infty$
with exponent $\rho$.

In the following list we state those Tauberian results needed in the
proofs of the preparatory lemmas and propositions of this sections,
as well as for the main theorems. 

\begin{itemize}
\item [(KTT)]\textbf{Karamata's~Tauberian~Theorem} (\cite{Feller:71},
\cite{Seneta:76}) Let $\left(b_{n}\right)_{n\geq0}$ be a non-negative
sequence such that for all $s>0,\; B\left(s\right):=\sum_{n\geq0}b_{n}e^{-ns}<\infty.$
Let $L$ be slowly varying at $\infty,$ and $\rho\in[0,\infty).$
Then \[
B\left(s\right)\sim\left(\frac{1}{s}\right)^{\rho}L\left(\frac{1}{s}\right)\qquad\mathrm{as}\; s\searrow0,\]
if and only\[
\sum_{k=0}^{n-1}b_{k}\sim\frac{1}{\Gamma\left(\rho+1\right)}n^{\rho}L\left(n\right)\qquad\mathrm{as\;}n\to\infty.\]
If $\left(b_{n}\right)$ is eventually monotone and $\rho>0,$ then
both are equivalent to \[
b_{n}\sim\frac{1}{\Gamma\left(\rho\right)}n^{\rho-1}L\left(n\right)\qquad\mathrm{as\;}n\to\infty.\]

\item [(KL)] \textbf{Karamata's Lemma} (\cite{Feller:71,Karamata33b}).
If $\left(a_{n}\right)$ is a regularly varying sequence with exponent
$\rho$ and if $p\geq-\rho-1,$ then \[
\lim_{n\to\infty}\frac{n^{p+1}a_{n}}{\sum_{k\leq n}k^{p}a_{k}}=p+\rho+1.\]

\item [(UA)] \textbf{Uniformly asymptotic} (\cite{Seneta:76}) Let $\left(p_{n}\right)$
and $\left(q_{n}\right)$ be two positive sequences with $p_{n}\to\infty$
and $\frac{p_{n}}{q_{n}}\in\left[1/k,k\right]$, $k>0$. Then we have
for $L$ a slowly varying function\[
\lim_{n\to\infty}\frac{L\left(p_{n}\right)}{L\left(q_{n}\right)}=1.\]

\item [(EL)] \textbf{Erickson Lemma} (\cite{Erickson:70}) Let $L\nearrow\infty$
be a monotone increasing continuous slowly varying function. Let $a_{t}\left(x\right)$
be defined by $a_{t}\left(x\right):=L^{-1}\left(xL\left(t\right)\right)$
with $x\in\left(0,1\right),$ where $L^{-1}\left(\cdot\right)$ denoting
the inverse function of $L\left(\cdot\right)$. Then we have for every
fixed $x\in\left(0,1\right)$\[
a_{t}\left(x\right)=o\left(t\right)\quad\textrm{and}\quad a_{t}\left(x\right)\longrightarrow\infty\qquad\left(t\to\infty\right).\]

\end{itemize}

\section{Proof of main results\label{sec:Proof-of-main}}

\proc{Proof of Proposition \ref{neue pro}.} Let $A\in\mathcal{A},\;0<\mu\left(A\right)<\infty,$
be a uniformly returning set, and let the functions $U\left(s\right),\; Q\left(s\right),\; s>0,$
be defined as Laplace transforms\begin{eqnarray}
Q\left(s\right) & := & \sum_{n=0}^{\infty}\frac{\mu\left(A\cap\left\{ \varphi>n\right\} \right)}{\mu\left(A\right)}e^{-ns}\label{Q1gleichung}\\
U\left(s\right) & := & \sum_{n=0}^{\infty}\nu\left(T^{-n}A\right)e^{-ns},\nonumber \end{eqnarray}
where $\nu$ denotes the probability measure with density $f\in\mathcal{P_{\mu}}.$

Since \[
A_{n}:=\bigcup_{k=0}^{n}T^{-k}A=\bigcup_{k=0}^{n}T^{-k}\left(A\cap\left\{ \varphi>n-k\right\} \right),\]
and the sets $T^{-k}\left(A\cap\left\{ \varphi>n-k\right\} \right),\;0\leq k\leq n,$
are disjoint, we have\[
\nu\left(A_{n}\right)=\int_{A}\sum_{k=0}^{n}\hat{T}^{k}\left(f\right)\cdot1_{A\cap\left\{ \varphi>n-k\right\} }\; d\mu,\quad n\geq0.\]
Thus,\[
\sum_{n=0}^{\infty}\nu\left(A_{n}\right)e^{-ns}=\int_{A}\left(\sum_{n=0}^{\infty}\hat{T}^{n}\left(f\right)e^{-ns}\right)\left(\sum_{n=0}^{\infty}1_{A\cap\left\{ \varphi>n\right\} }e^{-ns}\right)\; d\mu.\]
From\[
\hat{T}^{n}\left(f\right)\sim\frac{\nu\left(T^{-n}A\right)}{\mu\left(A\right)}\quad\textrm{as}\; n\to\infty\quad\mu\textrm{-a.s.}\quad\textrm{uniformly\; on}\; A,\]
it follows that\[
\sum_{n=0}^{\infty}\hat{T}^{n}\left(f\right)e^{-ns}\sim\frac{U\left(s\right)}{\mu\left(A\right)}\quad\textrm{as}\; s\to0\quad\mu\textrm{-a.s.}\quad\textrm{uniformly\; on}\; A.\]
This implies\[
\sum_{n=0}^{\infty}\nu\left(A_{n}\right)e^{-ns}\sim U\left(s\right)Q\left(s\right)\quad\textrm{as}\; s\to0.\]
Hence, since $\lim_{n\to\infty}\nu\left(A_{n}\right)=1,$ we obtain\begin{equation}
\frac{1}{s}\sim U\left(s\right)Q\left(s\right).\label{eq U}\end{equation}
If $b_{n}\sim n^{\beta}L\left(n\right)$ for $\beta\in\left[0,1\right)$,
$L$ denoting some slowly varying function, then due to (KL) we have\[
\sum_{k=0}^{n-1}\nu\left(T^{-k}A\right)\sim\frac{n^{1-\beta}}{\left(1-\beta\right)L\left(n\right)}\mu\left(A\right).\]
Thus, by (KTT) we obtain\[
Q\left(s\right)\sim\frac{1}{\Gamma\left(1-\beta\right)\mu\left(A\right)}\left(\frac{1}{s}\right)^{\beta}L\left(\frac{1}{s}\right).\]
Hence,\[
W_{n}\sim\frac{1}{\Gamma\left(1-\beta\right)\Gamma\left(1+\beta\right)}b_{n}.\]
Now let $W_{n}\sim n^{\beta}\tilde{L}\left(n\right)$ for $\beta\in\left[0,1\right)$,
$\tilde{L}$ denoting some slowly varying function. From (\ref{eq U}),
it follows by (KTT) \[
\sum_{k=0}^{n-1}\nu\left(T^{-k}A\right)\sim\frac{n^{1-\beta}}{\Gamma\left(2-\beta\right)\Gamma\left(1+\beta\right)\tilde{L}\left(n\right)}\mu\left(A\right).\]
Hence since \[
\sum_{k=0}^{n-1}\nu\left(T^{-k}A\right)\sim\mu(A)\sum_{k=0}^{n-1}\frac{1}{b_{k}},\]
$\left(b_{n}\right)$ is monotone, and $1-\beta>0$, we obtain by
(KTT)\[
\frac{1}{b_{n}}\sim\frac{n^{-\beta}}{\Gamma\left(1-\beta\right)\Gamma\left(1+\beta\right)\tilde{L}\left(n\right)}.\]
Thus,

\[
b_{n}\sim\Gamma\left(1-\beta\right)\Gamma\left(1+\beta\right)W_{n}.\]
From this the assertion follows.\ep\medbreak

\proc{Proof of Proposition \ref{propo0}.} 

\textbf{ad (1)} It is known that for every regularly varying function
with exponent $\beta\in\mathbb{R}$ there exists a slowly varying
function such that $F\left(x\right)=x^{\beta}L\left(x\right)$ for
all $x>0.$ Therefore to prove the result in Proposition \ref{propo0},
it suffices to show\[
\frac{L\left(Y_{n}\right)}{L\left(n\right)}\stackrel{\mathbb{P}}{\longrightarrow}1.\]
 We have for all $\delta>0$ and $K>0$ with $\delta<K$\[
\liminf\mathbb{P}\left(\delta\leq\frac{Y_{n}}{n}\leq K\right)\geq1-C_{\delta,K}.\]
 Due to the uniform convergence theorem for slowly varying functions
(cf. \cite{Seneta:76}) we have, for all $\varepsilon>0$ there exists
$n_{0}:=n_{0}\left(\varepsilon\right)$ such that\[
n\geq n_{0}\qquad\Rightarrow\qquad\left|\frac{L\left(\lambda n\right)}{L\left(n\right)}-1\right|<\varepsilon\quad\textrm{for all }\lambda\in\left[\delta,K\right].\]
Hence, for sufficiently large $n,$\[
\mathbb{P}\left(\left|\frac{L\left(Y_{n}\right)}{L\left(n\right)}-1\right|\geq\varepsilon\right)\leq1-\mathbb{P}\left(\delta\leq\frac{Y_{n}}{n}\leq K\right).\]
This implies \[
\limsup\mathbb{P}\left(\left|\frac{L\left(Y_{n}\right)}{L\left(n\right)}-1\right|\geq\varepsilon\right)\leq C_{\delta,K}.\]
 Since $C_{\delta,K}\longrightarrow0$ as $\delta\to0$ and $K\to\infty$,
the first part of the proposition follows.

\textbf{ad (2)} Now let $Y=0$ and $\beta>0$. Without loss of generality
we assume that $F$ is a positive and locally bounded function on
$[0,\infty)$. Then from \cite[p. 28]{BinghamGoldieTeugels:89} we
know that \[
F'(x):=\inf\left\{ y\geq0:F\left(y\right)>x\right\} ,\qquad x\in[0,\infty),\]
defines an asymptotic inverse of $F$, i.e. 

\begin{itemize}
\item [(1)]$F'\left(F\left(x\right)\right)\sim x$ as $x\to\infty$,
\item [(2)] $F'$ is regular varying with exponent $1/\beta$. 
\end{itemize}
Since $F(x)>y$ implies $y\geq F'\left(x\right)$we have for $\epsilon>0$
\[
\frac{F\left(Y_{n}\right)}{F\left(n\right)}>\epsilon\implies Y_{n}\geq F'\left(\epsilon F\left(n\right)\right).\]
Hence for all $n\geq1$ we have\[
\mathbb{P}\left(\frac{F\left(Y_{n}\right)}{F\left(n\right)}>\epsilon\right)\leq\mathbb{P}\left(\frac{Y_{n}}{n}>\epsilon_{n}\right)\]
with $\epsilon_{n}:=F'\left(\epsilon F\left(n\right)\right)n^{-1}$.
Since $F\left(n\right)\to\infty$, $n\to\infty$, we have by the properties
(1) and (2) of the asymptotic inverse that $\lim_{n\to\infty}\epsilon_{n}=\epsilon^{1/\beta}$.
Thus for $n$ sufficiently large, we have\[
\mathbb{P}\left(\frac{F\left(Y_{n}\right)}{F\left(n\right)}>\epsilon\right)\leq\mathbb{P}\left(\frac{Y_{n}}{n}>\frac{\epsilon^{1/\beta}}{2}\right)\]
proving the second part of the proposition for $\beta>0$. The case
$\beta<0$ is reduced to the above case by considering $1/F$ instead
of $F$. 

\ep\medbreak

If $T$ is a nonsingular ergodic transformation on $\left(X,\mathcal{A,\mu}\right),$
the compactness theorem in \cite[Sec. 3.6]{Aaronson:97} implies that
$R_{n}\circ T-R_{n}\stackrel{\mu}{\longrightarrow}0$ and $R_{n}\stackrel{\nu}{\Longrightarrow}R$
for some $\nu\in\mathcal{P_{\mu}}$ is sufficient for $R_{n}\stackrel{\mathcal{L\left(\mu\right)}}{\Longrightarrow}R$.
Hence, before proving the main theorems we state and prove the following
two lemmata.

\begin{lem}
\label{lem3}Let $A\in\mathcal{A}$ be a set of positive finite measure
$\mu\left(A\right)$ and $F\left(t\right)\to\infty$, $t\to\infty$,
be a regularly varying function with exponent $\beta\geq0.$ Then
we have\begin{equation}
\frac{1}{F\left(n\right)}\left(F\left(Z_{n}\circ T\right)-F\left(Z_{n}\right)\right)\;\stackrel{\mu}{\longrightarrow}\;0.\label{eq:101}\end{equation}

\end{lem}
\proc{Proof.}  Due to the Representation Theorem of regularly varying
functions (c.f. \cite{Seneta:76}) we have for some $B>0$ \begin{equation}
F\left(x\right)=x^{\beta}\psi\left(x\right)\exp\left(\int_{B}^{x}\frac{\zeta\left(t\right)}{t}\; dt\right)\quad\textrm{for\; all}\; x\geq B,\label{eq:100}\end{equation}
where $\psi$ is a positive measurable function on $\left[B,\infty\right)$
and $\zeta$ a continuous function on $\left[B,\infty\right)$ such
that\[
\psi\left(x\right)\longrightarrow C\in\left(0,\infty\right)\quad\textrm{und}\quad\zeta\left(x\right)\longrightarrow0\qquad\left(x\to\infty\right).\]
 Without loss of generality we assume that there exists $\delta\in\left(0,1\right)$
such that \[
\left|\zeta\left(t\right)\right|<\delta\quad\textrm{for\; all}\; t\geq B.\]
By (\ref{eq:100}) we have \[
F\left(x\right)\sim Cx^{\beta}\exp\left(\int_{B}^{x}\frac{\zeta\left(t\right)}{t}\; dt\right)\qquad\left(x\to\infty\right).\]
We define $\widetilde{F}\left(x\right):=Cx^{\beta}\exp\left(\int_{B}^{x}\frac{\zeta\left(t\right)}{t}\; dt\right)$
for $x\geq B$ and show first that the claim in the lemma holds for
$\widetilde{F}$. 

Let $\varphi$ be the first return time to the set $A$ defined as
in (\ref{returntime}). For $\varepsilon>0$ we define\[
K_{\varepsilon,n}:=\left\{ \varphi\leq n\;\wedge\;\frac{1}{\widetilde{F}\left(n\right)}\left|\widetilde{F}\left(Z_{n}\circ T\right)-\widetilde{F}\left(Z_{n}\right)\right|\geq\varepsilon\right\} \quad\left(n\in\mathbb{N}\right).\]
Since\begin{equation}
Z_{n}\left(T\left(x\right)\right)=\left\{ \begin{array}{ll}
Z_{n}\left(x\right)-1, & x\in\left\{ \varphi\leq n\right\} \cap T^{-\left(n+1\right)}A^{c},\\
n, & x\in T^{-\left(n+1\right)}A,\end{array}\right.\label{z(n)-z(n-1)}\end{equation}
we conclude\begin{eqnarray*}
K_{\varepsilon,n} & \subset & \left(\left\{ \varphi\leq n\right\} \cap T^{-\left(n+1\right)}A^{c}\cap\left\{ \frac{1}{\widetilde{F}\left(n\right)}\left(\widetilde{F}\left(Z_{n}\right)-\widetilde{F}\left(Z_{n}-1\right)\right)\geq\varepsilon\right\} \right)\\
 &  & \cup\left(T^{-\left(n+1\right)}A\cap\left\{ \frac{1}{\widetilde{F}\left(n\right)}\left(\widetilde{F}\left(n\right)-\widetilde{F}\left(Z_{n}\right)\right)\geq\varepsilon\right\} \right).\end{eqnarray*}
Note that $\widetilde{F}$ is a monotone increasing regularly varying
function and that $Z_{n}\to\infty$ $\mu$--a.s.  Therefore ,\[
\frac{1}{\widetilde{F}\left(n\right)}\left(\widetilde{F}\left(Z_{n}\right)-\widetilde{F}\left(Z_{n}-1\right)\right)\leq\frac{1}{\widetilde{F}\left(Z_{n}\right)}\left(\widetilde{F}\left(Z_{n}\right)-\widetilde{F}\left(Z_{n}-1\right)\right)\longrightarrow0\;\mu\textrm{-a.s.}\]
This implies for all $\nu\in\mathcal{P}_{\mu}$\begin{equation}
\lim_{n\to\infty}\nu\left(\left\{ \varphi\leq n\right\} \cap T^{-\left(n+1\right)}A^{c}\cap\left\{ \frac{1}{\widetilde{F}\left(n\right)}\left(\widetilde{F}\left(Z_{n}\right)-\widetilde{F}\left(Z_{n}-1\right)\right)\geq\varepsilon\right\} \right)=0.\label{gleichung1}\end{equation}
 For $n\geq B$ large enough let $\Omega_{n}:=\left\{ \omega:Z_{n}\left(\omega\right)\geq B\right\} $.
Then for all $\omega\in\Omega_{n}$ \\
${\displaystyle \widetilde{F}\left(n\right)-\widetilde{F}\left(Z_{n}\left(\omega\right)\right)=C\exp\left(\int_{B}^{Z_{n}\left(\omega\right)}\frac{\zeta\left(t\right)}{t}\; dt\right)\times}$\[
\;\quad\quad\quad\quad\quad\;\quad\left[n^{\beta}\exp\left(\int_{Z_{n}\left(\omega\right)}^{n}\frac{\zeta\left(t\right)}{t}\; dt\right)-\left(Z_{n}\left(\omega\right)\right)^{\beta}\right].\]
Since $\left|\zeta\left(t\right)\right|<\delta$ on $\left[B,\infty\right)$,
there exists a constant $C_{\delta}$, such that \[
\widetilde{F}\left(n\right)-\widetilde{F}\left(Z_{n}\left(\omega\right)\right)\leq C_{\delta}\left(n^{\beta+\delta}-\left(Z_{n}\left(\omega\right)\right)^{\beta+\delta}\right)=:E.\]
Now we distinguish two cases\textbf{.}\\
\textbf{Case 1}: For $\beta+\delta\geq1$, by the Mean-Value Theorem,
we have \[
E\leq C_{\delta}\left(\beta+\delta\right)n^{\beta+\delta-1}\left(n-Z_{n}\left(\omega\right)\right).\]
\textbf{Case 2}: For $\beta+\delta<1$ we have $E\leq C_{\delta}\left(n-Z_{n}\left(\omega\right)\right).$
Hence,\begin{eqnarray*}
T^{-\left(n+1\right)}A\cap\left\{ \frac{\widetilde{F}\left(n\right)-\widetilde{F}\left(Z_{n}\right)}{\widetilde{F}\left(n\right)}\geq\varepsilon\right\}  & \cap & \Omega_{n}\\
 & \subset & T^{-\left(n+1\right)}A\cap\left\{ n-Z_{n}\geq c_{n}\varepsilon\right\} \\
 & = & \bigcup_{c_{n}\varepsilon\leq k\leq n-1}\left\{ Z_{n}=n-k\right\} \cap T^{-\left(n+1\right)}A\\
 & = & \bigcup_{c_{n}\varepsilon+1\leq k\leq n}T^{-\left(n-k+1\right)}\left(A\cap\left\{ \varphi=k\right\} \right),\end{eqnarray*}
where in Case 1 we set $c_{n}:=\frac{\widetilde{F}\left(n\right)}{n^{\beta+\delta-1}C_{\delta}\left(\beta+\delta\right)}=\frac{n^{1-\delta}L\left(n\right)}{C_{\delta}\left(\beta+\delta\right)}$
for some slowly varying function $L$ and in Case 2 $c_{n}:=C_{\delta}\widetilde{F}\left(n\right).$
By the choice of $\delta$, in both cases we have $c_{n}\to\infty$
as $n\to\infty.$ Using the invariance of $\mu$ we obtain \begin{eqnarray*}
\mu\bigg(T^{-\left(n+1\right)}A\cap\left\{ \frac{1}{\widetilde{F}\left(n\right)}\left(\widetilde{F}\left(n\right)-\widetilde{F}\left(Z_{n}\right)\right)\geq\varepsilon\right\}  & \cap & \Omega_{n}\bigg)\\
 & \leq & \sum_{k=\left\lfloor c_{n}\varepsilon+1\right\rfloor }^{n}\mu\left(A\cap\left\{ \varphi=k\right\} \right)\\
 & \le & \mu\left(A\cap\left\{ \varphi\geq\left\lfloor c_{n}\varepsilon+1\right\rfloor \right\} \right)\\
 &  & \longrightarrow0\quad\textrm{für}\quad n\to\infty.\end{eqnarray*}
This gives for $\nu\in\mathcal{P}_{\mu}$\[
\lim_{n\to\infty}\nu\left(T^{-\left(n+1\right)}A\cap\left\{ \frac{1}{\widetilde{F}\left(n\right)}\left(\widetilde{F}\left(n\right)-\widetilde{F}\left(Z_{n}\right)\right)\geq\varepsilon\right\} \cap\Omega_{n}\right)=0.\]
Since $\lim_{n\to\infty}\nu\left(\Omega_{n}^{c}\right)=0$ we conclude\[
\lim_{n\to\infty}\nu\left(T^{-\left(n+1\right)}A\cap\left\{ \frac{1}{\widetilde{F}\left(n\right)}\left(\widetilde{F}\left(n\right)-\widetilde{F}\left(Z_{n}\right)\right)\geq\varepsilon\right\} \right)=0.\]
Using this, the fact $\lim_{n\to\infty}\nu\left(\left\{ \varphi>n\right\} \right)=0$,
and equation (\ref{gleichung1}) we finally conclude for all $\nu\in\mathcal{P}_{\mu}$
\begin{eqnarray*}
\lim_{n\to\infty}\nu\left(\left\{ \frac{\left|\widetilde{F}\left(Z_{n}\circ T\right)-\widetilde{F}\left(Z_{n}\right)\right|}{\widetilde{F}\left(n\right)}\geq\varepsilon\right\} \right) & \leq & \lim_{n\to\infty}\nu\left(K_{\varepsilon,n}\right)+\lim_{n\to\infty}\nu\left(\left\{ \varphi>n\right\} \right)\\
 & = & 0.\end{eqnarray*}
This gives (\ref{eq:101}) for $\widetilde{F}$.

We are left to show that (\ref{eq:101}) holds for arbitrary $F$.
For \[
Y_{n}:=F\left(Z_{n}\circ T\right)-F\left(Z_{n}\right)\quad\textrm{and}\quad\widetilde{Y}_{n}:=\widetilde{F}\left(Z_{n}\circ T\right)-\widetilde{F}\left(Z_{n}\right)\]
we first show\[
\frac{Y_{n}}{\widetilde{F}(n)}-\frac{\widetilde{Y}_{n}}{\widetilde{F}(n)}\stackrel{\mu}{\longrightarrow}0.\]
 Indeed, since $F\left(x\right)\sim\widetilde{F}\left(x\right)$ for
$x\to\infty$ and $Z_{n}\to\infty$ $\mu\textrm{-a.s.}$ we find for
a.e. $\omega\in X$ and $\varepsilon>0$ a number $n_{0}:=n_{0}\left(\omega,\varepsilon\right)$
such that\[
\left(1-\frac{\varepsilon}{2}\right)\widetilde{F}\left(Z_{n}\left(\omega\right)\right)\leq F\left(Z_{n}\left(\omega\right)\right)\leq\left(1+\frac{\varepsilon}{2}\right)\widetilde{F}\left(Z_{n}\left(\omega\right)\right)\quad\forall\; n\geq n_{0}\]
and\[
\left(1-\frac{\varepsilon}{2}\right)\widetilde{F}\left(Z_{n}\circ T\left(\omega\right)\right)\leq F\left(Z_{n}\circ T\left(\omega\right)\right)\leq\left(1+\frac{\varepsilon}{2}\right)\widetilde{F}\left(Z_{n}\circ T\left(\omega\right)\right)\quad\forall\; n\geq n_{0}.\]
This implies for all $n\geq n_{0}$ \[
\left|\widetilde{Y}_{n}\left(\omega\right)-Y_{n}\left(\omega\right)\right|\leq\frac{\varepsilon}{2}\left(\widetilde{F}\left(Z_{n}\circ T\left(\omega\right)\right)+\widetilde{F}\left(Z_{n}\left(\omega\right)\right)\right)\]
and by the monotonicity of $\widetilde{F}$ we have\[
\left|\frac{\widetilde{Y}_{n}\left(\omega\right)}{\widetilde{F}\left(n\right)}-\frac{Y_{n}\left(\omega\right)}{\widetilde{F}\left(n\right)}\right|\leq\varepsilon.\]
This shows that $\left(\frac{\widetilde{Y}_{n}\left(\omega\right)}{\widetilde{F}\left(n\right)}-\frac{Y_{n}\left(\omega\right)}{\widetilde{F}\left(n\right)}\right)\to0$
$\mu$-a.e. and consequently this convergence holds also locally stochastic
with respect to $\mu$. Using $F\left(n\right)\sim\widetilde{F}\left(n\right)$
this gives (\ref{eq:101}) for $F$.

\ep\medbreak

\begin{lem}
\label{lem4}Let $A\in\mathcal{A}$ be a set of positive finite measure
$\mu\left(A\right),$ then\[
\Psi_{n}\circ T-\Psi_{n}\;\stackrel{\mu}{\longrightarrow}\;0.\]

\end{lem}
\proc{Proof.} Let $\varphi$ be the first return time to the set
$A$. Let $\varepsilon>0$ be given, and let \[
K_{\varepsilon,n}:=\left\{ \varphi\leq n\;\wedge\;\left|\Psi_{n}\circ T-\Psi_{n}\right|\geq\varepsilon\right\} .\]
 Choose $n$ large enough such that $\frac{\mu\left(A\right)}{W_{n}}<\varepsilon.$
By (\ref{z(n)-z(n-1)}) we have \begin{eqnarray*}
K_{\varepsilon,n} & \subset & T^{-\left(n+1\right)}A\;\cap\;\left\{ n-Z_{n}\geq\frac{W_{n}}{\mu\left(A\right)}\varepsilon\right\} \\
 & = & \bigcup_{\frac{W_{n}}{\mu\left(A\right)}\varepsilon\leq k\leq n-1}\left\{ Z_{n}=n-k\right\} \cap T^{-\left(n+1\right)}A\\
 & = & \bigcup_{\frac{W_{n}}{\mu\left(A\right)}\varepsilon+1\leq k\leq n}T^{-\left(n-k+1\right)}\left(A\cap\left\{ \varphi=k\right\} \right).\end{eqnarray*}
Using the invariance of $\mu$ we obtain\[
\mu\left(K_{\varepsilon,n}\right)\leq\sum_{k=\left[\frac{W_{n}}{\mu\left(A\right)}\varepsilon+1\right]}^{n}\mu\left(A\cap\left\{ \varphi=k\right\} \right)\leq\mu\left(A\cap\left\{ \varphi\geq\left[\frac{W_{n}}{\mu\left(A\right)}\varepsilon+1\right]\right\} \right),\]
and therefore \[
\lim_{n\to\infty}\mu\left(K_{\varepsilon,n}\right)=0.\]
 This implies \[
\lim_{n\to\infty}\nu\left(K_{\varepsilon,n}\right)=0\quad\mathrm{for\; all}\;\nu\in P_{\mu}.\]
Since also $\lim_{n\to\infty}\nu\left(\left\{ \varphi>n\right\} \right)=0$
we have\[
\lim_{n\to\infty}\nu\left(\left\{ \left|\Psi_{n}\circ T-\Psi_{n}\right|\geq\varepsilon\right\} \right)\leq\lim_{n\to\infty}\nu\left(K_{\varepsilon,n}\right)+\lim_{n\to\infty}\nu\left(\left\{ \varphi>n\right\} \right)=0.\]
\ep\medbreak

Now we are in the position to give the proof of Theorem \ref{theo1}.

\proc{Proof of Theorem \ref{theo1}.} Let $A$ be a uniform set for
some $f\in\mathcal{P}_{\mu}$. We suppose without loss of generality
that $L$ is strictly monoton increasing and continuous. Then for
every $x\in\left(0,1\right)$ we have\begin{eqnarray*}
\nu\left(\frac{L\left(Z_{n}\right)}{L\left(n\right)}\leq x\right) & = & \nu\left(Z_{n}\leq a_{n}\left(x\right)\right)\\
 & = & \int_{A}\sum_{k=0}^{\left\lfloor a_{n}\left(x\right)\right\rfloor }\1_{A\cap\left\{ \varphi>n-k\right\} }\hat{T}^{k}\left(f\right)\; d\mu,\end{eqnarray*}
where $\nu$ denotes the probability measure with density $f\in\mathcal{P_{\mu}}$
and $a_{n}\left(x\right):=L^{-1}\left(xL\left(n\right)\right)$. By
(EL) we obtain \[
a_{n}\left(x\right)\rightarrow\infty\quad\textrm{and}\quad\frac{a_{n}\left(x\right)}{n}\rightarrow0\qquad\textrm{for}\; n\to\infty.\]
 Using the monotonicity of the sequence $\left(\1_{A\cap\left\{ \varphi>n\right\} }\right)$
we obtain by the asymptotic in (\ref{eq:AsymptAaronson}) on the one
hand that\begin{eqnarray*}
\nu\left(\frac{L\left(Z_{n}\right)}{L\left(n\right)}\leq x\right) & \leq & \int_{A}\1_{A\cap\left\{ \varphi>n-\left\lfloor a_{n}\left(x\right)\right\rfloor \right\} }\sum_{k=0}^{\left\lfloor a_{n}\left(x\right)\right\rfloor }\hat{T}^{k}\left(f\right)\; d\mu\\
 & \sim & \mu\left(A\cap\left\{ \varphi>n-\left\lfloor a_{n}\left(x\right)\right\rfloor \right\} \right)\cdot L\left(\left\lfloor a_{n}\left(x\right)\right\rfloor \right).\end{eqnarray*}
This and (UA) imply\[
\limsup\nu\left(\frac{L\left(Z_{n}\right)}{L\left(n\right)}\leq x\right)\leq x.\]
On the other hand we derive in a similar way \begin{eqnarray*}
\nu\left(\frac{L\left(Z_{n}\right)}{L\left(n\right)}\leq x\right) & \geq & \int_{A}\1_{A\cap\left\{ \varphi>n\right\} }\sum_{k=0}^{\left\lfloor a_{n}\left(x\right)\right\rfloor }\hat{T}^{k}\left(f\right)\; d\mu\\
 & \sim & \mu\left(A\cap\left\{ \varphi>n\right\} \right)\cdot L\left(\left\lfloor a_{n}\left(x\right)\right\rfloor \right).\end{eqnarray*}
This gives the opposite inequality\[
\liminf\nu\left(\frac{L\left(Z_{n}\right)}{L\left(n\right)}\leq x\right)\geq x.\]
 Finally the theorem follows from Lemma \ref{lem3} (for the case
$\beta=0$) by the compactness theorem.\ep\medbreak

For the proof of Theorem \ref{theo2} we need the following lemma.

\begin{lem}
\label{return lemma}Let $A\in\mathcal{A}$ with $0<\mu\left(A\right)<\infty$
be a \emph{}uniformly returning set, $W_{n}\sim L(n)$, where $L$
satisfies the properties stated in \emph{(EL),} and $x\in\left(0,1\right)$
fixed. Then for all $\varepsilon\in\left(0,1\right)$ there exists
$n_{0}$ such that for all $n\geq n_{0}$ and $k\in\left[n-a_{n}\left(x\right),n\right]$
we have uniformly on $A$\[
\left(1-\varepsilon\right)\frac{1}{W_{n}}\leq\hat{T}^{k}\left(f\right)\leq\left(1+\varepsilon\right)^{2}\frac{1}{W_{n}},\]
where $a_{n}\left(x\right)$ is defined as in \emph{(EL)}.
\end{lem}
\proc{Proof.}  Due to Proposition \ref{neue pro} we have \[
W_{n}\hat{T}^{n}\left(f\right)\;\longrightarrow\;1\quad\mu-\textrm{a.e. \; uniformly\; on}\; A.\]
Thus, for all \emph{$\varepsilon\in\left(0,1\right)$} there exists
\emph{}$k_{0}:=k_{0}\left(\varepsilon\right)$ such that we have uniformly
on $A$\[
\left(1-\varepsilon\right)\frac{1}{W_{k}}\leq\hat{T}^{k}\left(f\right)\leq\left(1+\varepsilon\right)\frac{1}{W_{k}}\quad\textrm{for\; all}\; k\geq k_{0}.\]
 By (EL) there exist $n_{1}$ and $n_{2}$ such that \[
n-a_{n}\left(x\right)\geq k_{0}\quad\textrm{for \; all}\; n\geq n_{1}\quad\textrm{and}\quad\frac{1}{W_{\left[n-a_{n}\left(x\right)\right]}}\leq\left(1+\varepsilon\right)\frac{1}{W_{n}}\quad\textrm{for \; all}\; n\geq n_{2}.\]
Let us denote $n_{0}:=\max\left\{ n_{1},n_{2}\right\} $. Then by
monotonicity of $W_{n}$ we obtain uniformly on $A$ that, for all
$n\geq n_{0}$ and $k\in\left[n-a_{n}\left(x\right),n\right],$\[
\left(1-\varepsilon\right)\frac{1}{W_{n}}\leq\left(1-\varepsilon\right)\frac{1}{W_{k}}\leq\hat{T}^{k}\left(f\right)\leq\left(1+\varepsilon\right)\frac{1}{W_{k}}\leq\frac{\left(1+\varepsilon\right)}{W_{\left[n-a_{n}\left(x\right)\right]}}\leq\frac{\left(1+\varepsilon\right)^{2}}{W_{n}}.\]
\ep\medbreak

\proc{Proof of Theorem \ref{theo2}.} Let $W_{n}\sim L\left(n\right)$
as $n\to\infty,$ without loss of generality we may assume that $L$
is monotone increasing and continuous. We have for every fixed $x\in\left(0,1\right)$\begin{eqnarray*}
\nu\left(\frac{L\left(n-Z_{n}\right)}{L\left(n\right)}\leq x\right) & = & \nu\left(Z_{n}\geq n-a_{n}\left(x\right)\right)\\
 & = & \int_{A}\sum_{n-a_{n}\left(x\right)\leq k\leq n}1_{A\cap\left\{ \varphi>n-k\right\} }\hat{T}^{k}\left(f\right),\end{eqnarray*}
where $a_{n}\left(x\right)=L^{-1}\left(xL\left(n\right)\right).$
Let $\varepsilon>0$. From Lemma \ref{return lemma} it follows that,
for sufficiently large $n$\[
\nu\left(\frac{L\left(n-Z_{n}\right)}{L\left(n\right)}\leq x\right)\leq\left(1+\varepsilon\right)^{2}\frac{1}{W_{n}}W_{\left[a_{n}\left(x\right)\right]}\sim\left(1+\varepsilon\right)^{2}x.\]
 Similarly for sufficiently large $n$,\[
x\left(1-\varepsilon\right)\sim\left(1-\varepsilon\right)\frac{1}{W_{n}}W_{\left[a_{n}\left(x\right)\right]}\leq\nu\left(\frac{L\left(n-Z_{n}\right)}{L\left(n\right)}\leq x\right).\]
Both inequalities give

\[
x\left(1-\varepsilon\right)\leq\liminf\nu\left(\frac{L\left(n-Z_{n}\right)}{L\left(n\right)}\leq x\right)\leq\limsup\nu\left(\frac{L\left(n-Z_{n}\right)}{L\left(n\right)}\leq x\right)\leq\left(1+\varepsilon\right)^{2}x.\]
Since $\varepsilon$ was arbitrary, we obtain \begin{equation}
\nu\left(\frac{L\left(n-Z_{n}\right)}{L\left(n\right)}\leq x\right)\longrightarrow x\quad\textrm{as}\quad n\to\infty\qquad\textrm{for all}\; x\in\left(0,1\right).\label{eq:prob0}\end{equation}
To show that the above convergence still holds if we replace $L$
by the wandering rate, we firstly point out that (\ref{eq:prob0})
in particular implies \begin{equation}
n-Z_{n}\longrightarrow\infty\;\textrm{ in probability (w.r.t.}\;\nu\textrm{)}\quad\textrm{as}\quad n\to\infty.\label{prob1}\end{equation}
 This can be seen as follows. At first note that (\ref{eq:prob0})
is equivalent to\begin{equation}
\nu\left(n-Z_{n}\leq a_{n}\left(x\right)\right)\longrightarrow x\quad\textrm{as}\quad n\to\infty\qquad\textrm{for all}\; x\in\left(0,1\right).\label{prob2}\end{equation}
 Now we suppose that (\ref{prob1}) fails. Then there exists $\varepsilon>0$,
a monotone increasing sequence $t_{n}\nearrow\infty$ and an integer
$N$ such that \[
\nu\left(t_{n}-Z_{t_{n}}\leq N\right)\geq\varepsilon\qquad\textrm{for\; all}\; n\in\mathbb{N}.\]
For arbitrary but fixed $x\in\left(0,\varepsilon\right)$ we have
$\lim_{n\to\infty}a_{n}\left(x\right)=\infty$. Hence, there exists
$n_{0}\in\mathbb{N}$ with \[
a_{t_{n}}\left(x\right)\geq N\qquad\textrm{for all}\; n\geq n_{0}.\]
 Thus,\[
\nu\left(t_{n}-Z_{t_{n}}\leq a_{t_{n}}\left(x\right)\right)\geq\nu\left(t_{n}-Z_{t_{n}}\leq N\right)\geq\varepsilon\qquad\textrm{for\; all}\; n\geq n_{0}.\]
This implies\[
\lim_{n\to\infty}\nu\left(t_{n}-Z_{t_{n}}\leq a_{t_{n}}\left(x\right)\right)\geq\varepsilon,\]
 contradicting (\ref{prob2}). 

Finally, since $n-Z_{n}\to\infty$ in probability, it is clear that
the slowly varying function $L$ may be replaced by any function $L_{1}$with
$L_{1}\left(n\right)\sim C\cdot L\left(n\right),$ $C>0$, as $n\to\infty.$
From this and Lemma \ref{lem4} by the compactness result the theorem
follows.

\ep


\end{document}